\documentclass[10pt]{asme2ej}

\usepackage{epsfig} 

\usepackage{graphics}
\usepackage{color}
\usepackage{amsmath,amssymb}
\usepackage{bm}
\usepackage[caption=false,font=footnotesize]{subfig}
\usepackage{txfonts}
\newcommand{\figref}[1]{Fig.\,\ref{#1}}
\newcommand{\tabref}[1]{Tab.\,\ref{#1}}
\newcommand{\secref}[1]{Sec.\,\ref{#1}}
\newcommand{\eref}[1]{(\ref{#1})}
\newcommand{\aref}[1]{(\ref{#1})}
\newcommand{\Ker}{{\rm Ker}}

\title{
A Lumped-Parameter Model of Multiscale Dynamics in Steam Supply Systems
}

\author{Hikaru Hoshino
    \affiliation{
Department of Electrical Engineering\\
Kyoto University\\
Katsura, Nishikyo, Kyoto, 615-8510 Japan\\
Email: hoshino@dove.kuee.kyoto-u.ac.jp
    }	\\
\\
{\tensfb Yoshihiko Susuki\thanks{Y. Susuki is currently with Department of  Electrical and Information Systems, Osaka Prefecture University, 1-1 Gakuencho, Nakaku, Sakai, Osaka 599-8531, Japan and also with JST CREST, 4-1-8 Honcho, Kawaguchi, Saitama, 332-0012 Japan} }\\
Department of Electrical Engineering\\
Kyoto University\\
Katsura, Nishikyo, Kyoto, 615-8510 Japan\\
Email: {susuki@eis.osakafu-u.ac.jp}\\
\\
{\tensfb Takashi Hikihara}\\
Department of Electrical Engineering\\
Kyoto University\\
Katsura, Nishikyo, \\
Kyoto, 615-8510 Japan\\
Email: hikihara.takashi.2n@kyoto-u.ac.jp
}

\begin{document}

\maketitle    

\begin{abstract}
{\it 
This paper focuses on multiscale dynamics occurring in steam supply systems. 
The dynamics of interest are originally described by a distributed-parameter model for fast 
steam flows over a pipe network coupled with a lumped-parameter model for slow internal dynamics of boilers. 
We derive a lumped-parameter model for the dynamics through physically-relevant approximations. 
The derived model is then analyzed theoretically and numerically in terms of existence of normally hyperbolic invariant manifold in the phase space of the model. 
The existence of the manifold is a dynamical evidence that the derived model preserves the slow-fast dynamics, and  suggests a separation principle of short-term and long-term operations of steam supply systems, which is analogue to electric power systems. 
We also quantitatively verify the correctness of the derived model by comparison with brute-force simulation of the original model. 
}
\end{abstract}

%

\section{Introduction}
\label{sec:introduction}

In this paper we study a problem of mathematical modeling for dynamics occurring in steam supply systems 
that consist of distributed plants producing, consuming, and interchanging steam. 
This type of steam supply is crucial to realization of energy systems integration \cite{geidl07:_energ_hub_future,omalley13}, where multiple types of energy such as electricity, heat, and natural gas are managed consistently in order to satisfy specifications of stability, reliability, and energy efficiency. 
In this integration, the dynamics of energy transfer and conversion governed by different physical laws occur on a wide range of spatio-temporal 
scales \cite{geidl07:_energ_hub_future,omalley13}.   
Modeling such multi-scale dynamics is of basic importance for establishing the control principle of the integrated energy systems.  
The problem which we study in this paper originates from the interaction between electricity and heat supply systems. 
This paper is a substantially-enhanced version of our conference papers \cite{nolta2014,cdc15}. 
\footnote{In this paper, we newly present the unified 
description of multiscale dynamics of steam supply systems 
by proper scaling of governing equations in \secref{sec:physical_modeling}, derivation of inner-limit of the derived model in \secref{sec:dynamical_analysis}, and numerical simulations for phase-space analysis and for verification of the model in \secref{sec:numerical_example}, all of which are not reported in \cite{nolta2014,cdc15}.}

\begin{figure}[t]
\centering
 \includegraphics[width=0.6\hsize]{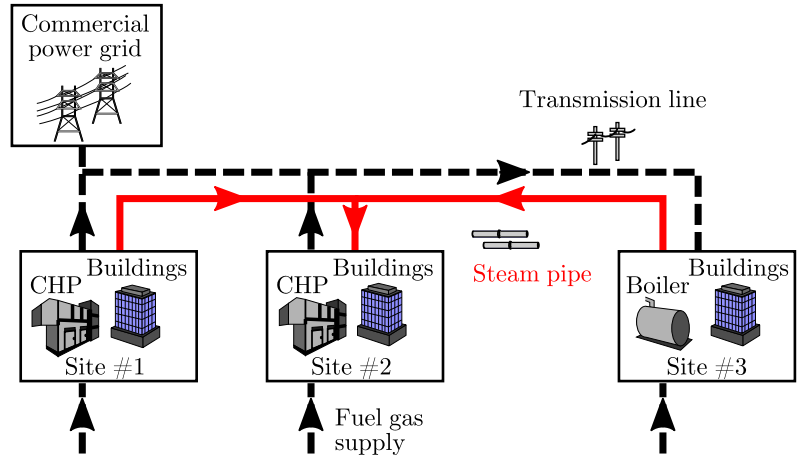}
\caption{Overview of steam supply systems with multiple CHP plants that we consider in this paper. It shows an example of three-site system. 
}%
\label{fig:overview}
\end{figure}

The interaction between electricity and heat appears in the spatial deployment of Combined Heat and Power (CHP) \cite{iea11} plants. 
The plants enable a local but spatially-distributed coupling between different energy systems because a CHP plant utilizes waste heat as a by-product of the conversion of fuel (natural gas, hydrogen, etc.) into electricity.
It is stated in \cite{geidl07:_energ_hub_future,omalley13} that this novel coupling would make it possible to design a new system architecture satisfying the above specifications. 
The overview of the target systems of this paper is shown in \figref{fig:overview} with an example of three-site system.   
In this figure, multiple boilers (including heat recovery boiler within CHP) are connected via a steam pipe network, and the produced heat is exchanged between different sites.  
This type of the steam supply systems has appeared in practice, e.g. district heating systems \cite{Rezaie12}. 
In such conventional steam-based systems, the primary objective of the CHP operation is to supply a desirable fixed amount of steam, and hence the operation does not necessarily contribute to the electricity supply \cite{Rezaie12}. 
However, in the view point of energy systems integration, it is possible to consider a novel operation of the CHP plants that contributes to both the steam and electricity supply by utilizing their ability of rapid electric response. 
In fact, it is proposed in \cite{iea11,shinji08,mueller14} that a CHP plant is used for compensating a variable output of renewable energy resources.  
This imposes a new problem on mathematical modeling for dynamics of steam supply systems 
against a large change of operating condition.

The contributions of this paper are twofold. 
The first one is to derive a \emph{lumped-parameter} model that captures stability and multiscale properties of steam supply. 
For internal dynamics in a single boiler or single plant, much work has been reported on lumped-parameter modeling \cite{astrom00,kim00,wen_ydstie09,bujak09}, whereas for dynamics of steam flows in pipes, partial differential equations or \emph{distributed-parameter} models \cite{osiadacz87,alobaid08,liu09,verma13} are normally used.
Although such models are crucial to plant design, detailed simulation, and commissioning, their simple  coupling is too complicated to reveal the system-wide dynamics of interest.   
They are originally described by the model for \emph{fast} steam flows over a pipe network coupled with the model for \emph{slow} internal dynamics of boilers.  
Namely, the lumped-parameter model of boilers is regarded as a slowly time-varying boundary condition of the {distributed}-{parameter} model of steam pipes. 
In this paper, through physically-relevant approximations, we newly derive a lumped-parameter model that contains multiscale dynamic characteristics of the steam pipes and boilers as well as a graph-theoretic property of the pipe network. 

The second contribution is to provide theoretical and numerical analyses of the derived model in terms of multiscale property of steam supply. 
The theoretical analysis is conducted with dynamical systems and graph theoretic methods \cite{kevorkian_cole96,wiggins94,nishiura02,iri69,pozrikidis14}. 
We obtain the inner limit of the derived model using the standard regular expansion method \cite{kevorkian_cole96} and locate a set of non-isolated equilibrium points of the inner-limit model. 
The set is thus proved to form a Normally Hyperbolic Invariant Manifold (NHIM) \cite{wiggins94,nishiura02} under mild technical conditions. 
The normal hyperbolicity characterizes the slow-fast vectorfield near the set, which is a dynamical evidence that the derived model preserves the slow-fast dynamics in the original model. 
Also, we conduct numerical simulations of the dynamics for an example of two-site system under a practical set of parameters.  
The correctness of the derived model is quantitatively verified by comparison with brute-force simulation of the original model, and the slow-fast dynamics near the NHIM are visualized.  
The existence of NHIM suggests a separation principle of short-term and long-term operations of steam supply systems, which is analogue to electricity supply operation \cite{machowski08}.

The rest of this paper is organized as follows. 
In \secref{sec:steam_supply_dynamics} we review 
the basic physical processes in steam supply systems. 
In \secref{sec:physical_modeling}, based on physical assumptions, we derive the lumped-parameter model. 
In \secref{sec:dynamical_analysis} we 
theoretically analyze the derived model and prove 
the existence of NHIM. 
In  
\secref{sec:numerical_example} we perform 
numerical simulations of 
the two-site system for providing a technological implication and for verifying the derived model. 
\secref{sec:conclusion} concludes this paper with a summary and future work.

\section{Physical processes in steam supply systems}
\label{sec:steam_supply_dynamics}

\begin{figure}[!t]
 \centering 
 \subfloat[Components of the two-site system]{\includegraphics[scale=0.25]{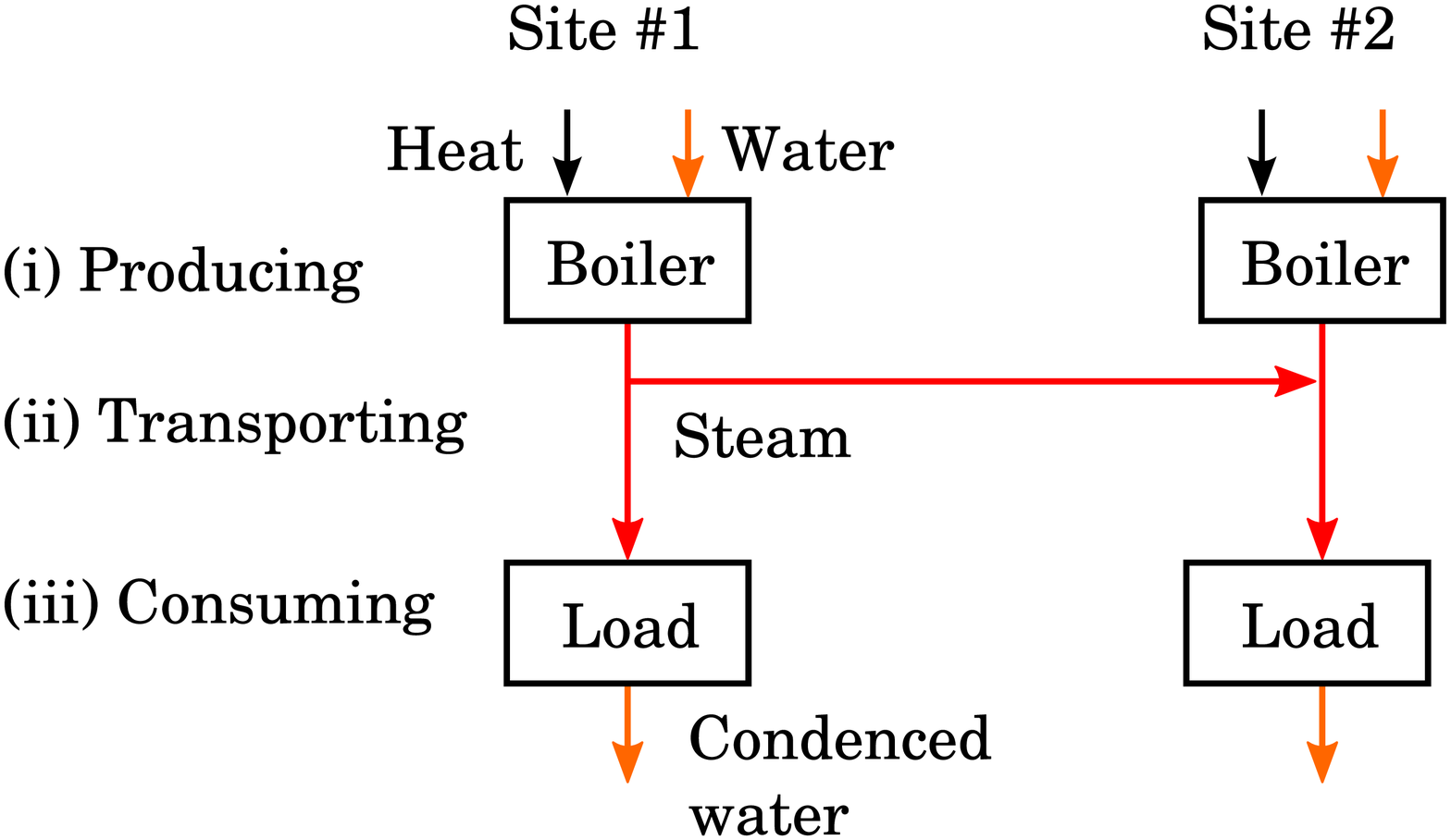}}
\hspace{1mm}
 \subfloat[Components of a boiler]{\includegraphics[scale=0.25]{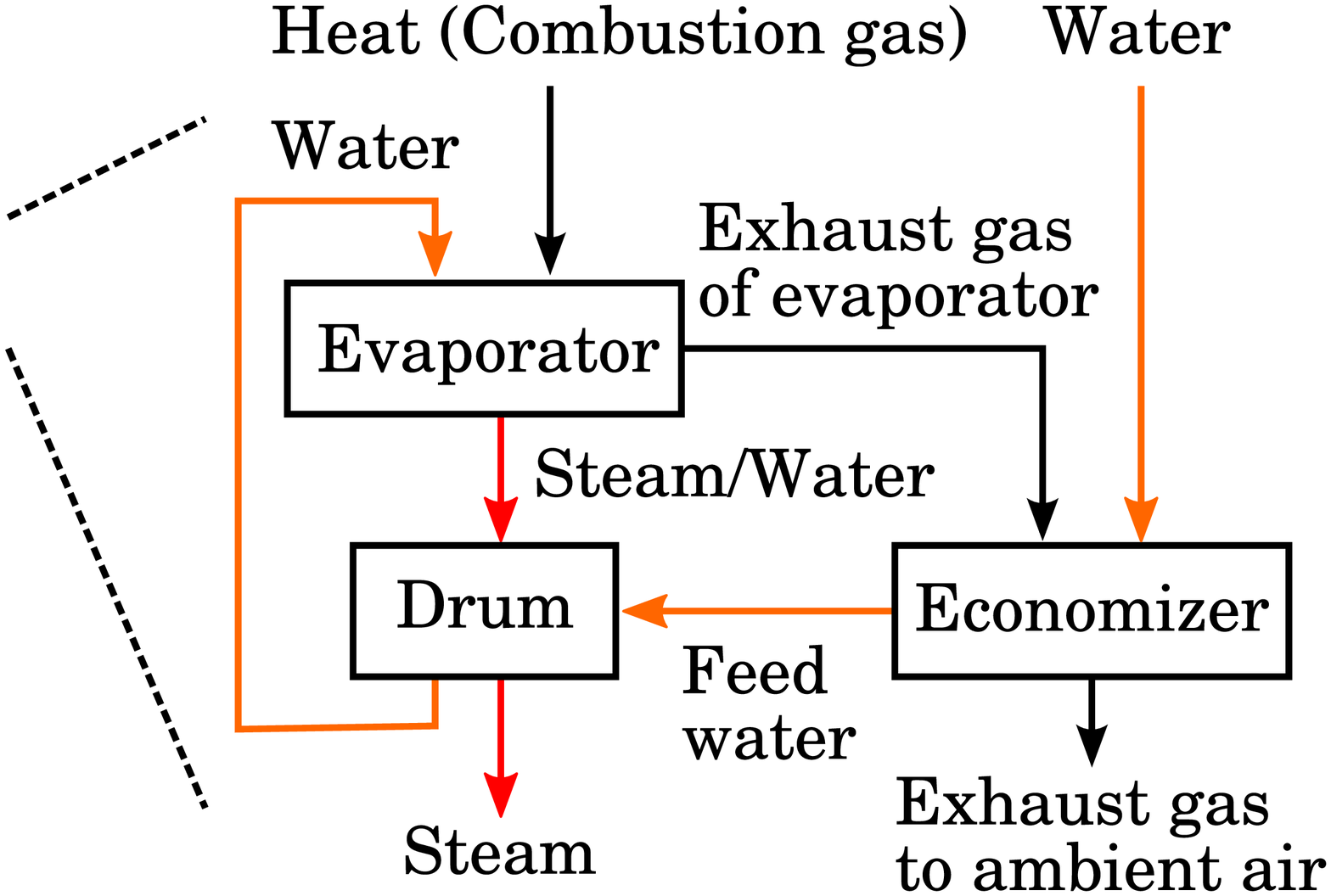}}
 \caption{Schematic diagram of the two-site steam supply. Each block in the figure shows (a) components of the two-site system and (b) components of a boiler. {\it Red}, {\it orange}, and {\it black} arrows describe the flow of steam, water, and heating gas, respectively.}
 \label{fig:heatsupply_schematic}
\end{figure}

This section describes the basic physical processes of steam supply systems that are related to (i) producing, (ii) transporting, and (iii) consuming steam. 
\figref{fig:heatsupply_schematic} illustrates the three processes for the two-site system. 
The blocks in the figure represent (a) components of the two-site system and (b) components of a boiler.  Each arrow describes the flow of steam, water, or heating gas. 
For (i), \figref{fig:heatsupply_schematic}\,(b) shows the production of steam from combustion gas and water in a boiler. 
Mixed steam and water are produced in the evaporator due to the boiling of water by high-temperature combustion gas.  
The mixed steam and water are supplied to the drum and stored at the phase equilibrium condition, while feedwater is supplied to the drum.
Thus, the saturated steam in the drum is brought to the outside of boiler.
For (ii), the transport of steam to a load in \figref{fig:heatsupply_schematic}\,(a) is realized due to the self-pressure of boiler.
The transport of steam between the multiple sites is also realized by controlling the difference of pressures between the two sites connected to a pipe. 
Finally, for (iii), the transported steam is consumed at each load.
The latent heat of the steam is extracted with a heat exchanger, and the resulting condensed water is returned to the boiler.

\section{Derivation of lumped-parameter model}
\label{sec:physical_modeling}

This section is devoted to the derivation of 
lumped-parameter model for dynamics of steam supply 
systems with the underlying physical processes in \secref{sec:steam_supply_dynamics}.  
The modeling procedure is exaplined in a dimensionless form in order to clearly describe the multiscale property of steam supply. 
The physical quantities with dimension are denoted by superscript $\ast$, and the reference quantities for scaling 
by subscript ${\rm r}$. 
The dimensionless time $t$ is scaled according to the time scale defined by the steam velocity $u_{\rm r}^\ast$ and the length scale $L_{\rm r}^\ast$ of the pipes: $t_{\rm r}^\ast:=L_{\rm r}^\ast/u_{\rm r}^\ast$. 
The detailed procedure of scaling and the reference quantities are presented in Appendix~\ref{sec:scaling}.

\subsection{Network description}
\label{sec:pn}
In this and next sections, we will consider a general steam supply system with arbitrary number of sites based on the graph theory \cite{iri69,pozrikidis14}. 
The topology of a steam supply system is described by a directed graph $\mathcal{G}=(\mathcal{V},\, \mathcal{L})$, where $\mathcal{V}$ stands for a finite set of vertices representing sites, and $\mathcal{L}$ for a finite set of links representing steam transporting pipes. 
For a link $l \in \mathcal{L}$, the tail (or head) vertex 
is denoted by $\partial^+l$ (or $\partial^-l$). 
For a vertex $v \in \mathcal{V}$, the set of outgoing (or incoming) links 
is denoted by  $\delta^+v$ (or $\delta^-v$). 
Below, the physical variables of boilers and pipes at each vertex and link are denoted by subscripts $v$ and $l$, respectively. 
The graph $\mathcal{G}$ is assumed to contain no self-loop and to be connected.
The assumption of no self-loop is relevant because a steam transporting pipe normally connects different two sites.  
Under this assumption, the incidence relation of a graph is completely represented 
with the incident matrix \cite{iri69}, and hence the matrix will be used in our modeling and analysis. 
The assumption of connected graph is intended for simplifying the presentation in this paper and does not lose generality of the modeling and analysis here.

\subsection{Physical assumptions}

As described in \secref{sec:steam_supply_dynamics}, the pressure and flow rate of steam are important physical quantities. 
In order to simply describe essential characteristics of the system-wide dynamics, we make the following assumptions:
\begin{enumerate}
 \renewcommand{\labelenumi}{(\theenumi)}
 \renewcommand{\theenumi}{A\arabic{enumi}}
 \item Temperatures of drum, evaporator and their wall are equal to temperature of mixed steam and water, i.e. the saturation temperature \cite{astrom00,kim00}.\label{as:boiler}
 \item Change of the volumes of steam and water is negligible when the water level of a drum is well regulated \cite{astrom00}.\label{as:water_level}
 \item Feedwater to a drum is at the condition of saturated liquid \cite{kim00}. \label{as:feedwater}
 \item Pressure drop in a pipe is evaluated by the Darcy-Weisbach equation for steady flow \cite{osiadacz87}. \label{as:darcy}
 \item No dominant effect of compressibility of steam appears on its velocity profile. This is relevant when the steam velocity is sufficiently smaller than the sound speed \cite{landau59}.\label{as:compressibility}
 \item No dominant effect of heat loss appears on the pressure drop and volumetric flow. This is relevant when the mass fraction of vapor in the fluid, i.e. the quality of steam is sufficiently close to one \cite{traviss71,kanro:english}. \label{as:condensation}
 \item Pressure drop in the site's components such as pressure regulators and valves are negligible.\label{as:local} 
\end{enumerate}
The validity of \aref{as:boiler} to \aref{as:feedwater} and \aref{as:darcy} have been tested in \cite{astrom00,kim00} and \cite{osiadacz87,alobaid08}, respectively. 
Thus, we mainly discuss the assumptions from \aref{as:compressibility} to \aref{as:local} in the rest of this paper.

\subsection{Steam boiler}
\label{sec:dynamics_boiler}

The dynamical model of a boiler is based on \cite{astrom00,kim00}. 
In the model, $V$ represents volume, $\rho$ density, $h$ specific enthalpy, $T$ temperature, and $m'$ mass flow rate. 
Furthermore, the three subscripts ${\rm s}$, ${\rm w}$, and ${\rm m}$ represent saturated steam, saturated water, and metal, respectively. 
The  total mass of metals of the drum and the evaporator is represented by  $m_{\rm t}$, and the specific heat of the metals by $C_{\rm p}$.
It is stated in \cite{astrom00} that the dynamics of pressure are well captured by global mass and energy balance. 
This is because the internal energy is rapidly released or absorbed due to the uniform boiling and condensation inside the drum and evaporator.
Thus, under the assumptions from \aref{as:boiler} to \aref{as:feedwater}, the dynamics of pressure $p_v$ at vertex $v \in \mathcal{V}$ are formulated as
\begin{equation}
 e_v(p_v) \frac{{\rm d} p_v}{{\rm d} t}= \epsilon_1 \left\{ Q_v' -m'_{{\rm s}v}h_{{\rm c}}(p_v) \right\}, \label{eq:boiler}
\end{equation}
where $Q'_v$ stands for the heat flow rate to the evaporator, and  $h_{\rm c}:=h_{\rm s}-h_{\rm w}$  corresponds to the enthalpy of condensation. 
The small parameter $\epsilon_1:={d_{\rm r}^\ast}^2L_{\rm r}^\ast/e_{\rm r}^\ast$ describes the slowness of the pressure dynamics \eref{eq:boiler} 
in terms of the time scale $t_{\rm r}^\ast$ of steam flow. 
From \cite{astrom00}, the coefficient $e_v(p_v)$ represents the rate of change of internal energy stored in the boiler against a change of pressure, 
given by 
\begin{align}
 e_v(p_v)= &h_{{\rm c}v}V_{{\rm s}v}\frac{\partial\rho_{{\rm s}v}}{\partial p_v} +\rho_{{\rm s}v}V_{{\rm s}v}\frac{\partial h_{{\rm s}}}{\partial p_v}  
 +\rho_{{\rm w}v}V_{{\rm w}v}\frac{\partial h_{{\rm w}v}}{\partial p_v} 
+m_{{\rm t}v}C_{\rm p}\frac{\partial T_{{\rm s}v}}{\partial p_v} -V_{{\rm s}v} -V_{{\rm w}v}.
  \label{eq:ei}
\end{align}
In this paper, according to \cite{astrom00}, the thermodynamic properties $h_{\rm s}$, $h_{\rm w}$, $\rho_{\rm s}$, $\rho_{\rm w}$, and $T_{\rm s}$ are evaluated from the steam table \cite{bergman11} and are represented as functions of pressure $p_v$, for example, $h_{{\rm s}v}=h_{\rm s}(p_v)$.

\subsection{Steam pipe}
\label{sec:dynamics_pipe}

\begin{figure}[t!]
 \centering
 \includegraphics[width=0.45\hsize]{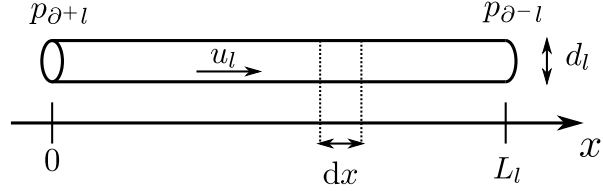}
 \caption{Schematic diagram of the steam transporting pipe $l$.}
 \label{fig:pipe}
\end{figure}

The transient steam flow in a pipe is described by the one-dimensional continuity equations of mass, momentum, and energy \cite{osiadacz87,alobaid08,liu09}. 
For each link $l \in \mathcal{L}$, the mass balance is given by 
\begin{equation}
 \dfrac{\partial\rho_l}{\partial t}+\dfrac{\partial}{\partial x}(\rho_l u_l)=0, \label{eq:mass_continuity}
\end{equation}
the momentum balance by
\begin{equation}
 \dfrac{\partial}{\partial t}(\rho_l u_l) + \dfrac{\partial}{\partial x}(\rho_l u_l^2)
 +\dfrac{1}{\epsilon_2}\dfrac{\partial p_l}{\partial x} +\lambda_l\dfrac{\rho_l u_l|u_l|}{2d_l} =0, \label{eq:momentum_continuity}
\end{equation}
and the energy balance by
\begin{equation}
 \label{eq:energy_continuity}
 \dfrac{\partial}{\partial t}(\rho_l h_l)+ \dfrac{\partial}{\partial x}(\rho_l h_l u_l)  = \dfrac{\partial p_l}{\partial t} +\epsilon_3 Q_{{\rm w}l},  
\end{equation}
where $u_l$ stands for the velocity of steam in a pipe $l$, and $x$ for the displacement variable along the pipe. 
The parameter $\lambda_l$ stands for the friction coefficient of the Darcy-Weisbach equation \cite{osiadacz87,kanro:english,bergman11} under \aref{as:darcy}. 
The parameters $d_l$ and $Q_{{\rm w}l}$ stand for the diameter of the pipe and the heat flow through walls, respectively. 
The parameters $\epsilon_2:=\rho_{\rm sr}^\ast{u_{\rm r}^\ast}^2/p_{\rm r}^\ast$ and $\epsilon_3:={d_{\rm r}^\ast}^2{L_{\rm r}^\ast}Q_{\rm wr}^\ast/Q_{\rm r}^{\prime\ast}$ in \eref{eq:momentum_continuity} and  \eref{eq:energy_continuity} are small, and they reflect the assumptions \aref{as:compressibility} and \aref{as:condensation}, respectively. 
Under the two assumptions, the above original equations are simplified through the incompressibility condition $\partial u/\partial x=0$. 
It is widely accepted that the low Mach number (described by $\epsilon_2$) implies an incompressible model \cite{landau59}. 
Further precise discussions are presented in \cite{muller98,principe09} for Navier-Stokes equations  and in \cite{felaco13} for one-dimensional flow equations.
The literature \cite{muller98,principe09,felaco13} shows that the simplification is relevant if both the constants $\epsilon_2$ and $\epsilon_3$ are sufficiently small.  
Hence, the energy equation \eref{eq:energy_continuity} is decoupled from the other equations, and the dynamics of steam flows are described by the momentum equation \eqref{eq:momentum_continuity} with the condition $\partial u/\partial x=0$: 
\begin{equation}
 \rho_l \dfrac{\partial u_l}{\partial t} +\dfrac{1}{\epsilon_2}\dfrac{\partial p_l}{\partial x}+ \dfrac{\lambda_l\rho_l u_l|u_l|}{2d_l}=0. 
\label{eq:dpdx}
\end{equation}
Therefore, by integrating \eqref{eq:dpdx} with respect to $x$ from $x=0$ to $x=L_l$ (see \figref{fig:pipe}), the following ordinary differential equation is derived: 
\begin{equation}
 L_l\rho_{{\rm av}l} \dfrac{{\rm d} u_l}{{\rm d} t}= \dfrac{p_{\partial^+l} -p_{\partial^- l}}{\epsilon_2} - \dfrac{\lambda_l \rho_{{\rm av}l} L_l u_l|u_l|}{2d_l}, \label{eq:velocity}
\end{equation}
where $\rho_{{\rm av}l}$ is given by
\begin{equation}
 \rho_{{\rm av}l}(t) := \dfrac{1}{L_l} \int_0^{L_l} \rho_l(x,t) {\rm d}x.
\end{equation}

\subsection{Site} 

At each site, equations~\eqref{eq:boiler} and \eqref{eq:velocity} are combined via the continuity equations of mass and energy. 
Under the assumption \aref{as:local} we have the following equations 
\cite{osiadacz87}: 
for each $v\in\mathcal{V}$ and all $t\in\mathbb{R}$, the mass balance is given by 
\begin{align}
  m'_{{\rm s}v}(t) = m'_{{\rm L}v}(t)
&+\sum_{l\in \delta^{+}v}\frac{\pi d_l^2}{4}\rho_l (0,t) u_l(t)  
  -\sum_{l\in \delta^-v}\frac{\pi d_l^2}{4}\rho_l (L_l,t) u_l(t) ,  \label{eq:massflow}
\end{align}
and the energy balance by 
\begin{align} 
 m'_{{\rm s}v}(t)h_{\rm s}(p_v(t)) 
 = m'_{{\rm L}v}(t)h_{{\rm L}v}(t) 
& +\sum_{l\in \delta^{+}v}\frac{\pi d_l^2}{4}\rho_l (0,t) u_l(t) h_l(0,t) \notag \\ 
& -\sum_{l\in \delta^-v}\frac{\pi d_l^2}{4}\rho_l(L_l,t) u_l(t) h_l(L_l,t) ,
  \label{eq:energyflow}
\end{align}
where $m'_{{\rm L}v}$ and  $h_{{\rm L}v}$ stand  
for 
the mass flow rate and specific enthalpy consumed by 
the load at vertex $v$.  
They are related to the consumption rate $Q'_{{\rm L}v}$ of heat as follows:
\begin{equation}
 m'_{{\rm L}v}h_{{\rm L}v} = Q'_{{\rm L}v} + m'_{{\rm L}v}h_{\rm w}(p_v). \label{eq:load}
\end{equation}
By multiplying both the sides of \eref{eq:massflow} by $h_{\rm w}(p_v)$ and using \eref{eq:energyflow}, we obtain
\begin{align}
 \label{eq:load_balance}
 m'_{{\rm s}v} h_{\rm c} (p_v) = Q'_{{\rm L}v} 
 &+\sum_{l\in \delta^{+}v}\frac{\pi d_l^2}{4}\rho_l (0,t)(h_l(0,t)-h_{\rm w}(p_v)) u_l \notag \\
 &-\sum_{l\in \delta^-v}\frac{\pi d_l^2}{4}\rho_l (L_l,t)(h_l(L_l,t)-h_{\rm w}(p_v))u_l.
\end{align}
The equation of $m'_{{\rm s}v} h_{\rm c}(p_v)$ is used in \eref{eq:boiler} and thereby combines the two equations \eqref{eq:boiler} and \eqref{eq:velocity}.

\subsection{Lumped-parameter model}
\label{sec:derived-model}

Finally, in order to derive a model in a self-consistent manner, it is necessary to determine the thermodynamic quantities in the pipes. 
While they are calculated in the original model by the equations \eqref{eq:mass_continuity} for mass and \eqref{eq:energy_continuity} for energy, 
we give their values by functions of pressure $p_v$ based on the assumptions \aref{as:condensation} and \aref{as:local}: 
\begin{align}
 \label{eq:thermodynamic_quantities}
 &h_l(0,t) := h_{\rm s}(p_{\partial^+l}(t)), \quad~ h_l(L_l,t) := h_{\rm s}(p_{\partial^-l}(t)), \notag \\ 
 &\rho_l(0,t) := \rho_{\rm s}(p_{\partial^+l}(t)), \quad~ \rho_l(L_l,t) := \rho_{\rm s}(p_{\partial^-l}(t)), \\
 &\rho_{{\rm av}l}(t) := \dfrac{\rho_{\rm s}(p_{\partial^+l}(t))+ \rho_{\rm s}(p_{\partial^-l}(t))}{2}. \notag
\end{align}
The relevance of these approximations will be discussed in Sec.\,\ref{ssec:comparison}. 

Consequently, the following model is derived for representing the steam supply dynamics: for each $v\in\mathcal{V}$ and $ l\in\mathcal{L}$,
\begin{subequations}
 \label{eq:network_model}
 \begin{align}
  & \frac{{\rm d} p_v}{{\rm d} t}= \dfrac{\epsilon}{e_v(p_v)} \left\{ Q'_v -Q'_{{\rm L}v} -\sum_{l\in \delta^{+}v}\frac{\pi d_l^2}{4}h_{\rm c}(p_v)\rho_{{\rm s}}(p_v) u_l +\sum_{l\in \delta^{-}v}\frac{\pi d_l^2}{4}h_{\rm c}(p_v) \rho_{{\rm s}}(p_v) u_l \right\}, \label{eq:network_model_boiler} \\
  & \frac{{\rm d} u_l}{{\rm d} t}= \dfrac{1}{\epsilon} \dfrac{2(p_{\partial^{+}l}-p_{\partial^{-}l})} {L_l\{\rho_{\rm s}(p_{\partial^+l})+\rho_{\rm s}(p_{\partial^-l})\}} -\frac{\lambda_l}{2d_l}u_l|u_l|, \label{eq:network_model_pipe}
 \end{align}
\end{subequations}%
where $\epsilon_1$ and $\epsilon_2$ were reset as $\epsilon$ by choosing the reference quantities as $e_{\rm r}^\ast:={d_{\rm r}^\ast}^2L_{\rm r}^\ast/\epsilon_2$.  
This resetting operation 
is relevant for a practical setting of parameters shown in \secref{sec:numerical_example}.

Here, we discuss the multiscale property of the derived model: see Sec.\,\ref{sec:dynamical_analysis} for its detailed analysis. 
The model \eqref{eq:network_model} 
includes a single small parameter $\epsilon$, 
and the parameters and functions on the right-hand sides are order of $1$. 
From \eqref{eq:network_model_boiler}, the pressure $p_v$ changes slowly in time due to the presence of the small parameter. 
In \eqref{eq:network_model_pipe}, the first and second terms on the right-hand side should be order of $1$. 
That is, the pressure difference $p_{\partial^+ l}-p_{\partial^- l}$ should be kept small (order of $\epsilon$). 
If this is not the case, for example if the pressure difference is $O(1)$, then the left-hand side of \eqref{eq:network_model_pipe} becomes $O(1/\epsilon)$, implying the steam velocity $u_l$ becomes large. 
This 
indicates the violation of \aref{as:compressibility}, and hence the derivation of the lumped-parameter model \eqref{eq:network_model} 
loses its validity.  
Thus, the smallness of $p_{\partial^+ l}-p_{\partial^- l}$ is necessary and will be assumed in the rest of this paper. 
Consequently, the short-term dynamics in $t\in [0, T]$ for $T=O(1)$ are 
described by the changes of $u_l$ in $O(1)$ and $p_v$ in $O(\epsilon)$. 
On the other hand, the long-term dynamics are related to the change of $p_v$ of $O(1)$, implying the changes of thermodynamic quantities such as $\rho_{\rm s}$.

\section{Theoretical analysis} 
\label{sec:dynamical_analysis}

This section shows that the derived model \eqref{eq:network_model} preserves the slow-fast dynamics in the original model. 
In phase-space geometric concepts, the presence of slow-fast dynamics can be described by the notion of normal hyperbolicity of invariant manifolds \cite{wiggins94,nishiura02}: slow dynamics along an invariant manifold and fast dynamics transversal to it. 
In order to theoretically analyze \eqref{eq:network_model}, here we simplify the model through the standard regular expansion method \cite{kevorkian_cole96}. 
As will be shown later, it corresponds to the inner-limit of the derived model.   
Thus, we prove the existence of NHIM for the inner-limit model \eqref{eq:graph-model}.
This indicates that the slow-fast dynamics are involved in the model \eref{eq:network_model}, because a NHIM persists under a perturbation of vectorfield \cite{wiggins94,nishiura02}. 
Indeed, in \secref{sec:numerical_example}, we will numerically confirm the existence of NHIM for the model \eref{eq:network_model}.

\subsection{Notation} \label{ssec:notation}
The notation frequently used in the rest of this paper is summarized below. 
The symbol $\mathsf{T}$ stands for the transpose operation of a vector or matrix.  
For a matrix $\mathsf{A}$, $\mathrm{Im}(\mathsf{A})$ represents the image space of linear mapping represented by $\mathsf{A}$, and
$\mathrm{Ker}(\mathsf{A})$ the kernel space of the linear mapping $\mathsf{A}$.  
The symbol $\mathrm{diag}(\bm{v})$ stands for the diagonal matrix made from a vector $\bm{v}$.   
For a vector-valued function ${\bm f}=(f_1,\ldots,f_m)^\top$, its Jacobian is denoted by $D{\bm f}$. 
The constant vector ${\bm 1}$ stands for all-one vector. 
For an Euclidean space $E$, $E^\bot$ stands for the orthogonal space of $E$.

\subsection{Derivation of inner-limit model}

First, we apply the regular expansion method to the derived model \eref{eq:network_model}.  
The 
regular expansion starts with assuming a solution of \eref{eq:network_model} in the following form: 
\begin{align}
 \label{eq:expansion}
 &p_v(t,\epsilon)=p_{v}^{(0)}(t)+\epsilon p_v^{(1)}(t) +O(\epsilon^2),\quad \notag \\
 &u_l(t,\epsilon)=u_{l}^{(0)}(t)+\epsilon u_l^{(1)}(t) +O(\epsilon^2). 
\end{align}
By substituting them into \eref{eq:network_model} 
and equating the coefficient of each power of $\epsilon$, 
a series of differential equations is obtained. 
From the leading-order terms, we obtain the following conditions: for all $v\in \mathcal{V}$ and $l \in \mathcal{L}$, 
\begin{equation}
 \dfrac{{\rm d} p_v^{(0)}}{{\rm d}t} =0,\quad p_{\partial^+ l}^{(0)}=p_{\partial^- l}^{(0)}.
\end{equation}
Thus, $p_v^{(0)}$ does not depend on both time and site, and is henceforth denoted by $p_0$. 
For the next-order of the series, we obtain the governing equations that contain the graph-theoretic property of the target steam supply system.  
For the graph $\mathcal{G}=(\mathcal{V},\, \mathcal{L})$, by labeling the vertices and links as $\mathcal{V}:=\{ v_1,\,\dots,\,v_n \}$ and  $\mathcal{L}:=\{ l_1,\,\dots,\,l_m \}$, the first-order term of pressure and the zeroth-order term of volumetric flow of steam are described by 
\begin{align}
 \bm{\psi}&:= \begin{bmatrix}p_1^{(1)},\,\dots,\,p_n^{(1)}\end{bmatrix}^\top, \quad 
 \bm{q} := \dfrac{\pi}{4}\begin{bmatrix} d_1^2 u_1^{(0)}, \,\dots,\, d_m^2 u_m^{(0)}\end{bmatrix}^\top.
\end{align}
The incidence matrix ${\sf R}=(R_{ij}) \in \mathbb{R}^{n\times m}$ of the graph is defined by
\begin{equation}
 R_{ij} :=  
\begin{cases}
	1 & {\rm if}~ v_i = \partial^+ l_j \neq  \partial^- l_j, \\
	-1 & {\rm else~if}~ v_i = \partial^- l_j \neq \partial^+ l_j,  \\
	0 & {\rm otherwise}.
\end{cases}
\end{equation}
With the notation, the dynamics of $(\bm{\psi},\bm{q})$ are described as follows: 
\begin{equation}
 \begin{bmatrix}
  {\sf G} &0  \\ 0& {\sf H}
 \end{bmatrix}
\dfrac{{\rm d}}{{\rm d}t}
 \begin{bmatrix}
 \bm{\psi} \\ \bm{q}
 \end{bmatrix} 
=
\begin{bmatrix}
 0 & -{\sf R}\\
 {\sf R}^\top & 0
\end{bmatrix}
 \begin{bmatrix}
  {\bm \psi} \\ {\bm q}
 \end{bmatrix} +
 \begin{bmatrix}
  {\bm s} \\ -{\bm f}({\bm q})
 \end{bmatrix},
 \label{eq:graph-model}
\end{equation}
with
\begin{align}
 {\sf G}&:= {\rm diag}\left( \dfrac{e_1(p_0)}{h_{\rm c}(p_0)\rho_{\rm s}(p_0)},\,\dots,\,\dfrac{e_n(p_0)}{h_{\rm c}(p_0)\rho_{\rm s}(p_0)}\right), \\
 {\sf H}&:= {\rm diag}\left( \dfrac{4\rho_{\rm s}(p_0)L_1}{\pi d_1^2},\,\dots,\,\dfrac{4\rho_{\rm s}(p_0)L_m}{\pi d_m^2} \right), \\
 {\bm s}&:= \left( \dfrac{Q'_1-Q'_{{\rm L}1}}{h_{\rm c}(p_0)\rho_{\rm s}(p_0)},\,\dots,\,\dfrac{Q'_n-Q'_{{\rm L}n}}{h_{\rm c}(p_0)\rho_{\rm s}(p_0)} \right)^\top, \label{eq:definition_s} \\
 {\bm f}({\bm q})&:= \left(\dfrac{8\lambda_1L_1\rho_{\rm s}(p_0) q_1|q_1|}{\pi^2 d_1^5}, \dots, \dfrac{8\lambda_m L_m \rho_{\rm s}(p_0) q_m|q_m|}{\pi^2 d_m^5} \right)^\top.  \label{eq:definition_f}
\end{align}
Note that \eref{eq:graph-model} possibly has a unbounded solution, and thus the expansion \eref{eq:expansion} is not uniformly valid for all time $t\in \mathbb{R}$. 
Namely, there exists a finite $T$ such that the expansion is valid at all $t \in [0, T]$, and thus the limitation $\epsilon \to 0$ implies the inner limit \cite{kevorkian_cole96}. 
In this sense, we call \eref{eq:graph-model} as the inner limit model of 
\eqref{eq:network_model}.

\subsection{Characterization of invariant manifold}
\label{ssec:equilibrium_analysis}

In this subsection, we locate the invariant manifold in the inner-limit model \eqref{eq:graph-model} 
using graph theory: see Appendix~\ref{sec:gt} for its summary. 
As will be shown below, the manifold is located as a set of non-isolated equilibrium points. 
An equilibrium point $({\bm \psi}^\ast,\,{\bm q}^\ast)$ of \eqref{eq:graph-model} satisfies the following condition:  for given $\bm s$, $\sf R$, and $\bm f$,
\begin{subequations}
\label{eq:AEep}
 \begin{align}
  &{\sf R}{\bm q}^\ast = {\bm s}, \label{eq:AEep-a}\\
  &{\sf R}^\top {\bm \psi}^\ast= {\bm f}({\bm q}^\ast). \label{eq:AEep-b}
 \end{align}
\end{subequations}%
The above conditions can be regarded as the combination of two linear equations defined by ${\sf R}$ and ${\sf R}^\top$. 
Thus, they are analyzed using the image and kernel spaces of ${\sf R}$ and ${\sf R}^\top$.
First, \eqref{eq:AEep-a} has a solution if $\bm{s}$ belongs to the image of ${\sf R}$:
\begin{equation}
  {\bm s} \in {\rm Im}({\sf R})= ({\rm Ker}({\sf R}^\top))^\bot. 
  \label{eq:s-condition}
\end{equation}
Since \eqref{eq:s-condition} implies that $\bm{s}$ is orthogonal to ${\bm 1}$ (see \eref{eq:kerRt} in Appendix~\ref{sec:gt}), the sum of all the elements $s_1+\cdots +s_n$ should be zero. 
By the definition 
in \eqref{eq:definition_s}, this condition is equivalent to 
\begin{equation}
 \label{eq:entire_heat_balance}
 Q'_1 + \cdots +Q'_n = Q'_{{\rm L}1} + \cdots +Q'_{{\rm L}n}.  
\end{equation}
This clearly indicates that the sum of all generation and consumption of steam at vertices is zero. 
Even if the above condition does not hold, the following analysis in this section is still applicable by introducing a new state variable $(\bm{\psi}', \bm{q})$ with the following time-varying transformation: 
\begin{equation}
 \bm{\psi}' = \bm{\psi} - \left( \dfrac{s_1+\dots+s_n}{G_{11}+\dots+G_{nn}} t \right)  \bm{1}.
\end{equation}

Thus, under the condition (\ref{eq:s-condition}), the set of all solutions of \eqref{eq:AEep-a}, denoted by $Q^\ast({\bm s})$, is described as follows:
\begin{equation}
 Q^\ast({\bm s}) = \big\{ \bm{q_0} + {\bm q} ~\big|~ {\bm q} \in {\rm Ker}({\sf R}) \big\}, \label{eq:y-candidate}
\end{equation}
where $\bm{q_0$} stands for one of the solutions of \eref{eq:AEep-a}.
Second, \eref{eq:AEep-b} has a solution if ${\bm f}({\bm q})$ 
belongs to the image of ${\sf R}^\top$:
\begin{equation}
 {\bm f}({\bm q}) \in {\rm Im}({\sf R}^\top).  \label{eq:y-condition}
\end{equation}
This implies that ${\bm f}({\bm q})$ belongs to the subspace of cutsets, which is known as the Kirchhoff's law for tensions \cite{iri69}. 
Under the condition (\ref{eq:y-condition}) and for ${\bm q}^\ast\in Q^\ast({\bm s})$, the set $\Psi^\ast({\bm q}^\ast)$ of all solutions of \eref{eq:AEep-b} is described as follows:
\begin{equation}
 \Psi^\ast({\bm q^\ast}) = \big\{ \bm{\psi_0}+{\bm \psi} ~\big|~ {\bm \psi}\in {\rm Ker({\sf R}^\top)}= c\,{\bm 1} ~,c \in\mathbb{R} \big\}. 
\label{eq:psi-zerospace}
\end{equation}
where $\bm{\psi_0}$ stands for one of the solutions of \eref{eq:AEep-b}. 
Consequently, under the condition
\begin{equation}
  {\bm s} \in ({\rm Ker}({\sf R}^\top))^\bot, \quad {\rm Im}({\sf R}^\top)\cap {\bm f}(Q^\ast({\bm s}))\neq  \emptyset, \label{eq:equilibrium-cond}
\end{equation}
the set of equilibrium points $({\bm \psi}^\ast, {\bm q}^\ast)$ is 
described as follows:
\begin{align}
 \Big\{ ({\bm \psi}^\ast, {\bm q}^\ast) \in \mathbb{R}^{n+m} \Big |~ &{\bm \psi}^\ast={\bm \psi}_0+ c\,{\bm 1}, ~c \in \mathbb{R},~  {\bm q}^\ast \in {\bm f}^{-1}({\rm Im({\sf R}^\top))}\cap Q^\ast({\bm s} )  \Big\},
 \label{eq:equilibrium_set}
\end{align}
where ${\bm f}^{-1} ({\rm Im}({\sf R}^\top)\cap {\bm f}(Q^\ast({\bm s})) = {\bm f}^{-1} ({\rm Im}({\sf R}^\top) ) \cap Q^\ast({\bm s}) $ holds because ${\bm f}$ is bijective from the definition in  \eqref{eq:definition_f}. 
This fact also indicates the uniqueness of ${\bm q}^\ast$ if it exists \cite{iri69}.
As a result, the set \eqref{eq:equilibrium_set} of equilibrium points becomes a one-dimensional curve (line) in the state space of the inner-limit model \eqref{eq:graph-model}. 
This curve forms an invariant manifold consisting of all the non-isolated equilibrium points, which we will denote by $\mathcal{I}$. 

In order to discuss the technological implication in \secref{sec:numerical_example}, the dynamics of the model \eqref{eq:network_model} near $\mathcal{I}$ are described below.
Since the invariant manifold herein is one-dimensional, the slow dynamics will be characterized by one-dimensional reduced system on it. 
A perturbation to solutions of \eqref{eq:graph-model} along  $\mathcal{I}$ represents the uniform change of the pressures, and thus the slow dynamics correspond to the responses of pressure level of the entire system.  
On the other hand, the dynamics transverse to $\mathcal{I}$ represent the change of steam velocities and pressure fluctuations near the manifold, and the fast dynamics correspond to the transport of steam between the different sites.

\subsection{Proof of normal hyperbolicity of the invariant manifold}

In this subsection, we prove that the located invariant manifold $\mathcal{I}$ in \eqref{eq:equilibrium_set} is normally hyperbolic under certain conditions.  
An invariant manifold is called normally hyperbolic if the expansion or contraction rate of vectors transverse to the manifold dominates that of vectors tangent to the manifold. 
For a precise formulation of normal hyperbolicity, see \cite{wiggins94,hirsch70}. 
For the present discussion, since $\mathcal{I}$ consists of non-isolated equilibrium points, it is characterized by eigenvalues associated with the linearization of the model \eqref{eq:graph-model} at each equilibrium point. 
By substituting ${\bm \psi}={\bm \psi}^\ast+\bm{\Delta \psi}$ and  ${\bm q}={\bm q}^\ast+\bm{\Delta q}$ into \eqref{eq:graph-model}, the linearized system around $({\bm \psi}^\ast,{\bm q}^\ast)$ is obtained as follows: 
\begin{equation}
\dfrac{\rm d}{{\rm d}t}
\begin{bmatrix}
 \bm{\Delta \psi}\\ \bm{\Delta q}
\end{bmatrix} 
 = {\sf A}
\begin{bmatrix}
 \bm{\Delta\psi}\\ \bm{\Delta q}
\end{bmatrix},
\end{equation}
with 
\begin{equation}
{\sf A}:=
\begin{bmatrix}
 0 & -{\sf G}^{-1}{\sf R}\\ {\sf H}^{-1}\,{\sf R}^\top & -{\sf H}^{-1}D{\bm f}(\bm{q}^\ast)
\end{bmatrix}. 
\label{eq:linearized}
\end{equation}
In below, we will show that the center subspace of the linearized system \eref{eq:linearized} is one-dimensional and is tangent to $\mathcal{I}$. 
To do this, we analyze the eigenvector associated with zero eigenvalue under the following two assumptions. 
The first one is that the matrix $D{\bm f}(\bm{q}^\ast)$ is non-singular. 
This is relevant  if $q_i^\ast\neq 0$ for all $i=1,\dots, m$.  
Under the assumption, the 
eigenvector $(\bm{\Delta \psi_0}, \bm{\Delta q_0})$ 
associated with zero eigenvalue 
satisfies the following equations:
\begin{subequations}
 \label{eq:center_cond}
 \begin{align}
  &{\sf G^{-1}}{\sf R} \bm{\Delta q_0}={\bm 0} \label{eq:center-cond1},\\
  &\bm{\Delta q_0} = D{\bm f}(\bm{q}^\ast)^{-1} \,{\sf R}^\top \bm{\Delta \psi_0}. 
\label{eq:center-cond2} \end{align}
\end{subequations}
The second assumption states non-existence of pure imaginary eigenvalues of ${\sf A}$. 
As stated in \cite{brayton_moser64}, non-singularity of $D{\bm f}(\bm{q}^\ast)$ is closely related to the non-oscillating condition, that is, non-existence of pure imaginary eigenvalues. 
However, in this paper, we will simply make both the assumptions. 
Under the two assumptions, one can verify that the center subspace is spanned by eigenvector associated to zero eigenvalue and is tangent to $\mathcal{I}$. 
By substituting \eqref{eq:center-cond2} into \eqref{eq:center-cond1}, the following condition holds for $\bm{\Delta \psi_0}$: 
\begin{equation}
 {\sf R}{\sf \Sigma}{\sf R}^\top \bm{\Delta \psi_0} ={\bm 0}, \quad {\sf \Sigma}:= D{\bm f}(\bm{q}^\ast)^{-1}.
\end{equation}
Since ${\sf \Sigma}$ is diagonal by definition, 
${\sf R}{\sf \Sigma}{\sf R}^\top$ corresponds to the so-called Kirchhoff matrix \cite{pozrikidis14}. 
From the assumption of connected graph, the kernel of ${\sf R}{\sf \Sigma}{\sf R}^\top$ is represented by
\begin{equation}
 \{ c {\bm 1} \in \mathbb{R}^n ~ |~ c \in \mathbb{R} \}.
\end{equation}
From  \eqref{eq:center-cond2}, $\bm{\Delta q_0}={\bm 0}$ holds because this space is also the kernel of ${\sf R}^\top$, that is, ${\sf R}^\top \bm{\Delta \psi_0}={\bm 0}$. 
Note that, by discussion similar to above, one can verify that there is no generalized eigenvector associated with zero eigenvalue other than  $(\bm{\Delta \psi_0}, \bm{\Delta q_0})$. 
Thus, the center subspace of the linearized system \eref{eq:linearized} is explicitly represented  as follows:
\begin{equation}
 \Big\{ ( \bm{\Delta \psi},{\bm 0}) \in \mathbb{R}^{n+m} | \bm{\Delta \psi} = c\,{\bm 1}, ~c\in \mathbb{R} \Big\}.
\end{equation}
This clearly indicates that the center subspace is one-dimensional and is tangent to the invariant manifold $\mathcal{I}$. 
From the definition, if the two assumptions---non-singularity of $D{\bm f}(\bm{q}^\ast)$ and non-existence of pure imaginary eigenvalues of ${\sf A}$---are satisfied, then $\mathcal{I}$ is normally hyperbolic. 
The proof is thus completed. 

\begin{table*}[!t] %
 \centering
 \caption{%
 List of variables and parameters in the derived model \eqref{eq:network_model}.  
 The values used for numerical simulations in Sec.\,\ref{sec:numerical_example} are also presented. 
 }%
 \label{tab:parameters}
 \small
 \begin{tabular}{l c r c} \hline
 Meaning & Symbol & (Nominal) Value & Scaled value\\\hline
 Pressure of steam & $p$ & $800\,{\rm kPa}$ &1.0 (base)\\
 Density of saturated steam& $\rho_{\rm s}$ &$4.16\,{\rm kg/m^3}$ & 1.0 (base) \\
 Density of saturated water& $\rho_{\rm w}$ &$897\,{\rm kg/m^3}$  \\
 Specific enthalpy of saturated steam&$h_{\rm s}$ & $2768\,{\rm kJ/kg}$ & 14.3\\
 Specific enthalpy of saturated water&  $h_{\rm w}$ & $721\,{\rm kJ/kg}$ & 3.74\\
 Specific enthalpy of the feed water& $h_{\rm f}$ & \\
 Mass flow rate of saturated steam from drum &  $m'_{{\rm s}}$ & \\
 Mass flow rate of feed water to drum&  $m'_{{\rm f}}$ & \\
 Coefficient of pressure variation given in \eref{eq:ei} & $e$ & $3073\,{\rm J/Pa}$ & $1.8$ \\
 Temperature of saturated steam &  $T_{{\rm s}}$ & $443\,{\rm K}$ \\
 Total mass of evaporator and drum of a boiler& $m_{{\rm t}}$ & $50,000\,{\rm kg}$\\
 Specific heat of the metal of boiler & $C_{\rm p} $ & $0.4\,{\rm kJ/(K\cdot kg)}$ \\
 Total volume of steam & $V_{{\rm s}}$ & $10.0\,{\rm m^3}$ \\  
 Total volume of water &  $V_{{\rm w}}$ & $10.0\,{\rm m^3}$ \\  
 Velocity of steam in the pipe & $u$ & $30\,{\rm m/s}$ & 1.0 (base)\\
 Length of the steam pipe &  $L$& $200\,{\rm m}$ & 1.0 (base)\\
 Diameter of the steam pipe &  $d$& $0.2\,{\rm m}$ & 1.0 (base) \\
 Friction coefficient of the steam pipe&   $\lambda$& 0.016 & 16 \\
 Input rate of heat of a boiler &   $Q'$& \\
 Consumption rate of heat at a load&  $Q'_{{\rm L}}$& $5.0\,{\rm MJ/s}$ & 5.2 \\\hline
 \end{tabular}
\end{table*}

\section{Numerical simulations for two-site system}
\label{sec:numerical_example} 

This section demonstrates the slow-fast dynamics near the NHIM and verifies the correctness of the derived model \eref{eq:network_model} by numerical simulations for the two-site system. 
For the system, in \cite{scipaper14-en}, we studied dynamics of electricity supply based on the classical formulation of power system swing equations \cite{machowski08} by assuming that the steam supply system is ideally operated. 
In this paper, following \cite{scipaper14-en}, we discuss the technological implication for the operation of steam supply systems by comparing to the electricity supply operation \cite{machowski08}. 
For the minimal two-site system, the dimension of the derived model is three, and hence it is possible to perfectly visualize the state space of the model including a NHIM and to apply phase-space geometric concepts \cite{wiggins94} to it.

\begin{figure}[t!]
 \centering
 \subfloat[$e_i(p_i)$ ({\it solid} line) and the values of terms on the right-hand side of \eref{eq:ei}]{\includegraphics[width=0.75\hsize]{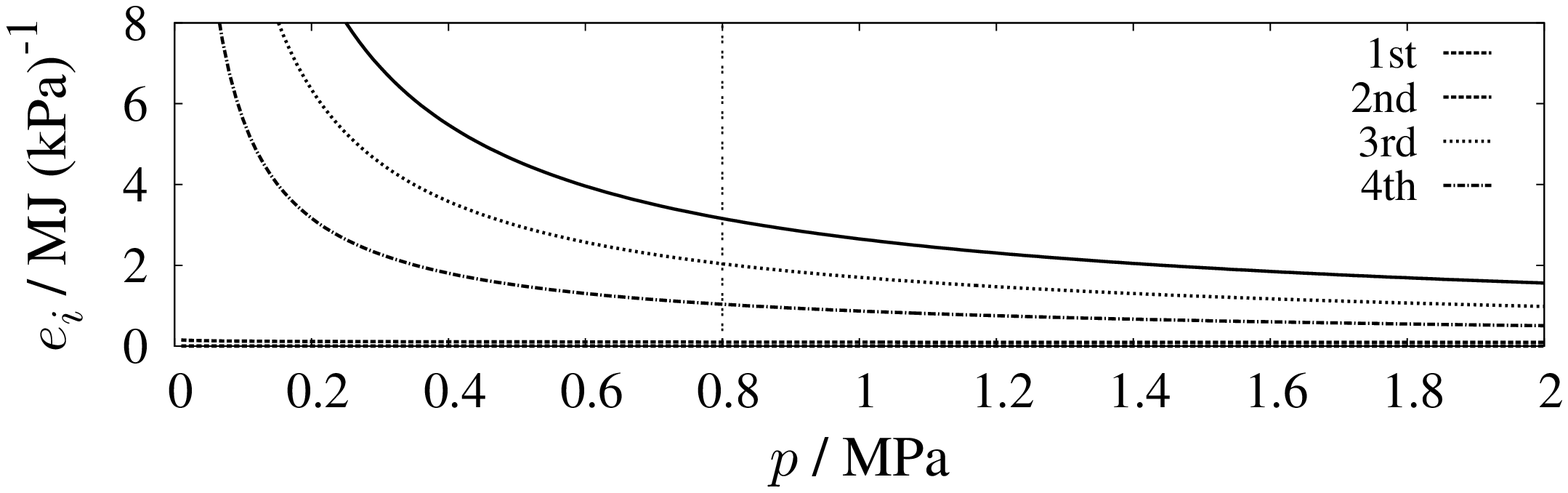}}\\
 \subfloat[$h_{\rm c}(p_i)\,\rho_{\rm s}(p_i)$]{\includegraphics[width=0.75\hsize]{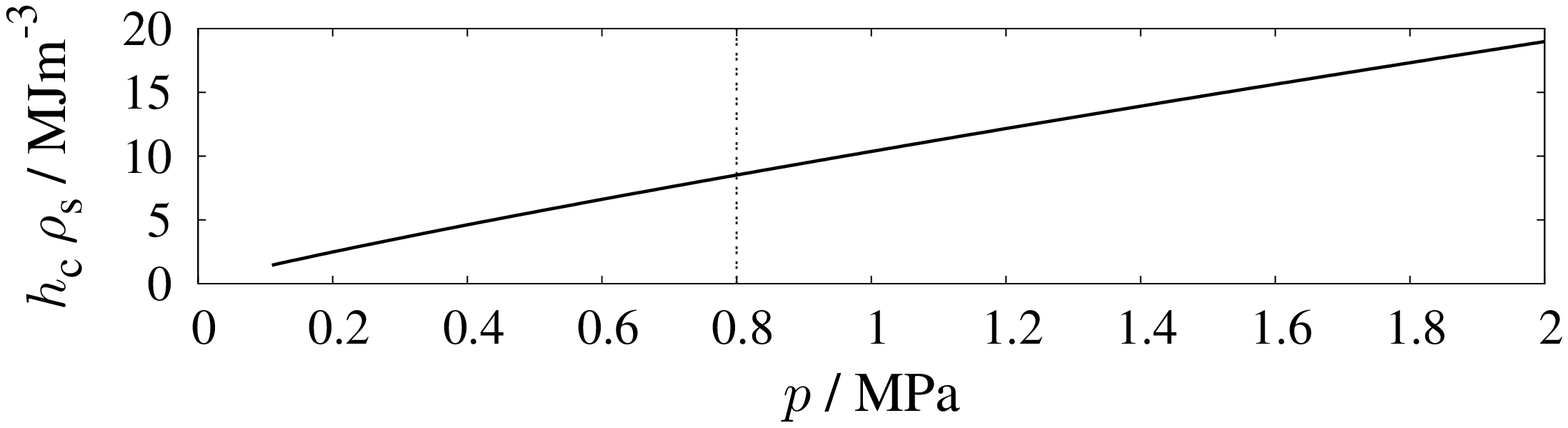}}
 \caption{Numerical values of (a) $e_i$ and (b) $h_{\rm c}\rho_{\rm s}$}
 \label{fig:functions}
\end{figure}

In \cite{astrom00} transient responses of a single boiler are examined experimentally as well as numerically under a setting of fuel profiles for a real plant.  
Following this, in 
this paper we provide responses of physical quantities under a step change and periodic change of the heat flow rates $Q'_i$ to boilers. 
Although the 
abrupt change of $Q'_i$ 
may not be possible in realistic operation, it provides basic information on the multiscale dynamics. 
The values of parameters used for the current simulations are shown in \tabref{tab:parameters} and based on district heating systems \cite{bujak09,kanro:english}.
According to \cite{wen_ydstie09}, the thermodynamic properties are calculated by using Xsteam package \cite{xsteam}. 
Numerical values of $e_i(p_i)$ and $h_{\rm c}(p_i)\rho_{\rm s}(p_i)$ in \eref{eq:network_model} are shown in \figref{fig:functions} for a practical range of pressure \cite{bujak09}: $0.03\,{\rm MPa} \le p_i \le 2\,{\rm MPa}$.  

\subsection{Time-response analysis} 
\label{subsec:time}

\begin{figure}[!t]
 \centering
 \includegraphics[width=0.75\hsize]{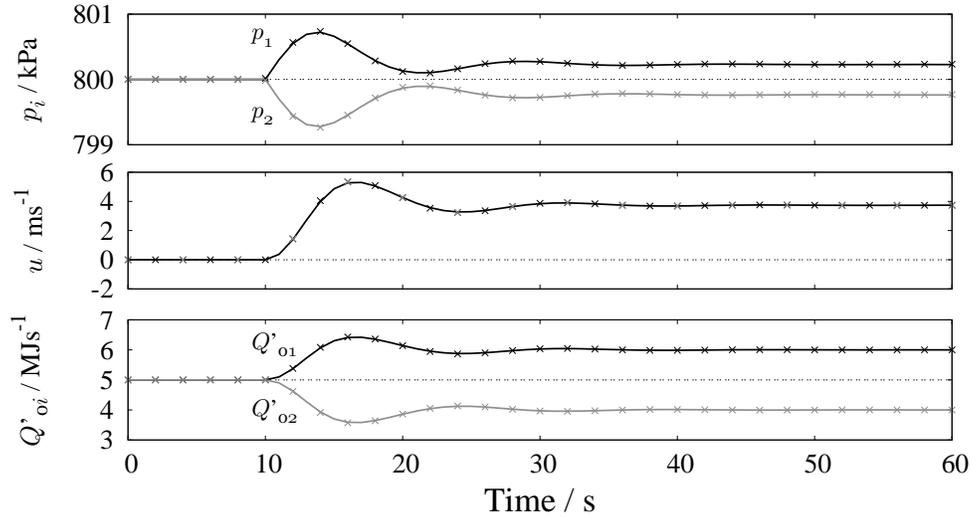}
 \caption{Responses of state variables $(p_1,\,p_2,\,u)$ and heat output rates $Q'_{{\rm o}i}:=m'_{{\rm s}i} h_{\rm c}(p_i)$ from the boilers. The \emph{solid} lines show the responses of  \eqref{eq:network_model}, and the \emph{points} ($\times$) show the sample plot of responses of \eqref{eq:graph-model}. They are initiated by a step change of the parameters from $(Q'_1,\,Q'_2)=(5\,{\rm MJ/s}, 5\,{\rm MJ/s})$ to $(Q'_1,\,Q'_2)=(6\,{\rm MJ/s}, 4\,{\rm MJ/s})$.}
 \label{fig:short-term}
\end{figure}

First, we simulate the short-term dynamics related to transport of steam between the two sites.  
Time-responses of physical quantities are provided under the following step change:
\begin{equation} 
 (Q'_1, Q'_2)= 
\begin{cases}
 (5\,{\rm MJ/s}, 5\,{\rm MJ/s}),&  t < 10\,{\rm s}  \\
 (6\,{\rm MJ/s}, 4\,{\rm MJ/s}), & t \ge 10\,{\rm s}  
\end{cases}
\label{eq:step_change}
\end{equation}
This implies an abrupt change of operating conditions of the boilers, and thereby $1\,{\rm MJ/s}$ surplus (or deficit) of heat is caused in site $\# 1$ (or site $\# 2$). 
Note that in both cases the condition \eqref{eq:entire_heat_balance} is satisfied. 
\figref{fig:short-term} shows step responses of the state variables $(p_1,\,p_2,\,u)$ of \eref{eq:network_model} and heat output rates $Q'_{{\rm o}i}$ from the boilers.
The system is initially at a steady operating condition with no transport of steam between the sites, and then $Q'_ i$ changes at $t=10\,{\rm s}$.
The pressures $p_1$ and $p_2$ and velocity $u$ move to a new operating condition after transients in a few tens of seconds.
At the new operating condition, the values of $p_1-p_2$ and $u$ become positive.
This clearly shows 
that the positive pressure drop $p_1-p_2$ induces the transport of steam from site $\# 1$ to site $\# 2$. 
The heat output rates $Q'_{\rm o1}$ and $Q'_{\rm o2}$ change symmetrically, and thus 
the surplus and deficit of heat at the two sites are compensated by the transport of steam.

\begin{figure}[!t]
 \centering
 \vspace{-12mm}
 \includegraphics[width=0.65\hsize]{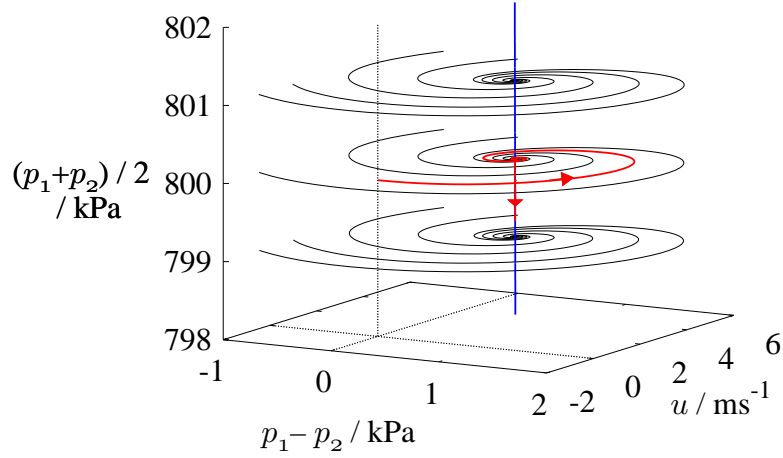} \\[-5mm]
 \caption{%
 Trajectories of the derived model \eqref{eq:network_model} with $(Q'_1,Q'_2)=(6\,{\rm MJ/s}, 4\,{\rm MJ/s})$. 
 The {\it red} trajectory corresponds to the time response presented in \figref{fig:short-term} with subsequent long-term response. The {\it blue} trajectory shows a invariant manifold located with direct numerical integration of \eref{eq:network_model}.
 }%
 \label{fig:3d-perturbed}
\end{figure}

\subsection{Phase-space analysis}

Second, we analyze the dynamics described by the model \eref{eq:network_model} from the viewpoint of phase space.
Especially, the long-term dynamics described by \eref{eq:network_model} are considered in this subsection. 
%
The model (\ref{eq:network_model}) for the two-site system has the three independent variables.
Based on the analysis in Sec.\,\ref{ssec:equilibrium_analysis}, we introduce the following variable transformation: 
\begin{equation}
\left[
\begin{array}{c}
(p_1+p_2)/2 \\
p_1-p_2 \\
u
\end{array}
\right]=
\left[
\begin{array}{ccc}
1/2 & 1/2 & 0\\
1 & -1 & 0 \\
0 & 0 & 1
\end{array}
\right]
\left[
\begin{array}{c}
p_1 \\
p_2 \\
u
\end{array}
\right].
\end{equation}
\figref{fig:3d-perturbed} shows trajectories of the model (\ref{eq:network_model}) under the parameter setting as $(Q'_1,\,Q'_2)=(6\,{\rm MJ/s}, 4\,{\rm MJ/s})$. 
The {\it red} trajectory corresponds to the time response presented in \figref{fig:short-term} with subsequent long-term response from $t=0\,{\rm s}$ to $10,000\,{\rm s}$. 
This trajectory shows typical \emph{slow-fast} dynamics as mentioned in 
\secref{sec:dynamical_analysis}. 
The mean pressure value $(p_1+p_2)/2$ does not change dominantly while the pressure difference $p_1-p_2$ and the velocity $u$ exhibit fast oscillations shown in \figref{fig:short-term}.
After the fast oscillation is settled, the mean pressure $(p_1+p_2)/2$ begins to decrease slowly. 
The {\it blue} trajectory in the figure shows a one-dimensional invariant manifold located with direct numerical integration of \eref{eq:network_model}. 
It is confirmed that the trajectories converge to the located invariant manifold while exhibiting the slow and fast dynamics as mentioned above.
This result becomes a numerical evidence that the located invariant manifold possesses 
the normal hyperbolicity.

\begin{figure}[!t]
 \centering
 \subfloat[Time responses]{\includegraphics[width=0.8\hsize]{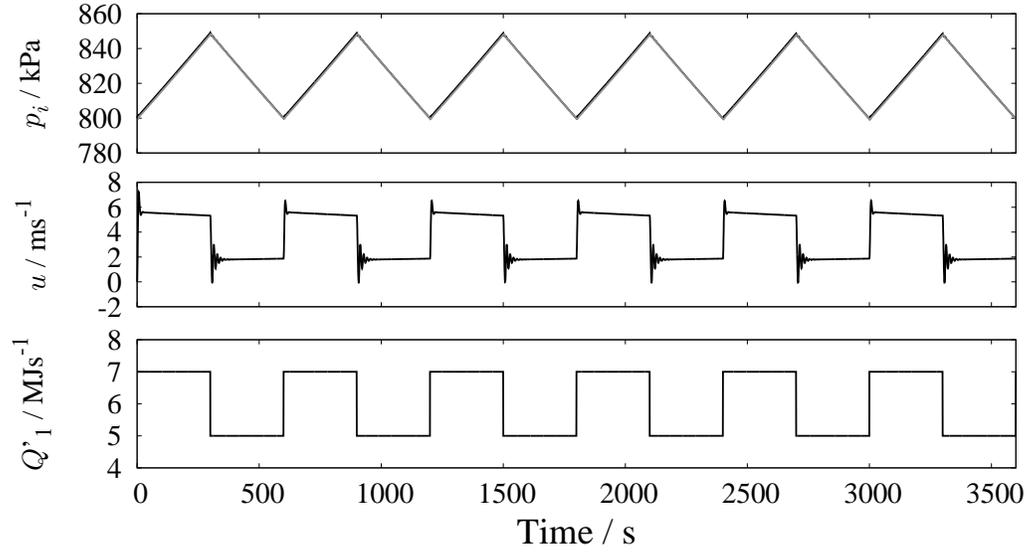}}\\[-12mm]
 \subfloat[Trajectory in the phase space]{\includegraphics[width=0.7\hsize]{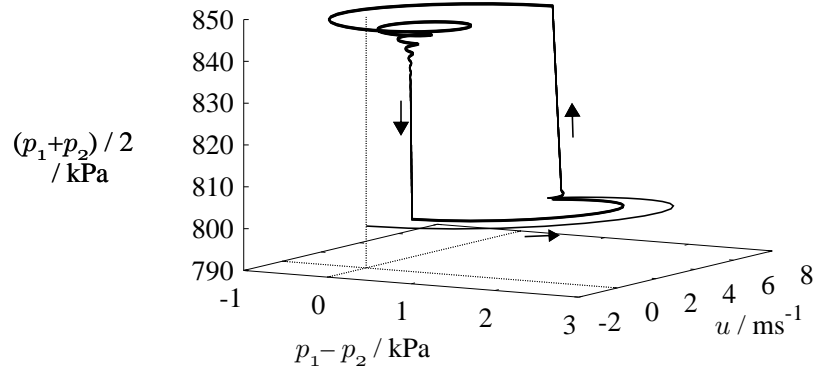}}
 \caption{%
 Long-term dynamics for the periodic change of the parameter $Q'_1$. 
 Numerical simulation of (a) time responses of the state variables and (b) 
 trajectory in the phase space are shown. 
 }%
 \label{fig:periodic}
\end{figure}

Also, we analyze the long-term dynamics related to the 
boilers' operation 
to offer several technological implications of the phase space analysis. 
Here, we consider the following periodic change of $Q'_{1}$: 
for $n=1,\,2,\dots$,
\begin{equation} 
 Q'_1= 
\begin{cases}
 7\,{\rm MJ/s}, & 600(n-1)\,{\rm s}\le t < 600n-300\,{\rm s},  \\
 5\,{\rm MJ/s}, & 600n-300\,{\rm s}\le t < 600n\,{\rm s},
\end{cases}
\end{equation}
and $Q'_2=4\,{\rm MJ/s}$. 
This periodic change is intended to the novel electricity-oriented operation of CHP plants mentioned in \secref{sec:introduction}, 
and the similar profiles of the fuel flow rate are shown in \cite{astrom00} as experiment data of a real plant. 
Under the above setting, the condition \eqref{eq:entire_heat_balance} is not satisfied.  
\figref{fig:periodic} shows (a) the time responses of the state variables and (b) the corresponding trajectory in the phase space. 
The appearing short-term and long-term dynamics are characterized by the NHIM. 
For the short-term regime, the fast motion towards to the NHIM guarantees that the pressure and velocity oscillation are settled after a change of operating condition.
This ensures that the amount of transport of steam between the two sites becomes bounded in the short-term regime.  
On the other hand, for 
the long-term regime, the mean pressure value changes due to the slow motion along the NHIM. 
Thus, the success of the novel electricity-oriented operation can be clarified as the existence of NHIM near which the separation of fast and slow motions holds.

It should be here noted that this type of slow-fast separation plays an important role in power system operation \cite{machowski08}.
Conventionally, the frequency dynamics in power systems can be classified into three stages with different time-scales and are regulated separately with different mechanisms \cite{machowski08}.
In this sense, we now identify a similar time-scale separation in the steam supply system, which has not been reported yet in literature. 
The time-scale separation enables independent operations of transport of steam between the sites in \emph{short-term} regime and \emph{long-term} supply-demand balancing of steam in order to maintain the pressure level in the system. 
The finding of the time-scale separation is thus expected to become a dynamical principle for operational design of steam supply.

\subsection{Comparison with the original model} 
\label{ssec:comparison}
\begin{figure}[!t]
 \centering
 \includegraphics[width=0.75\hsize]{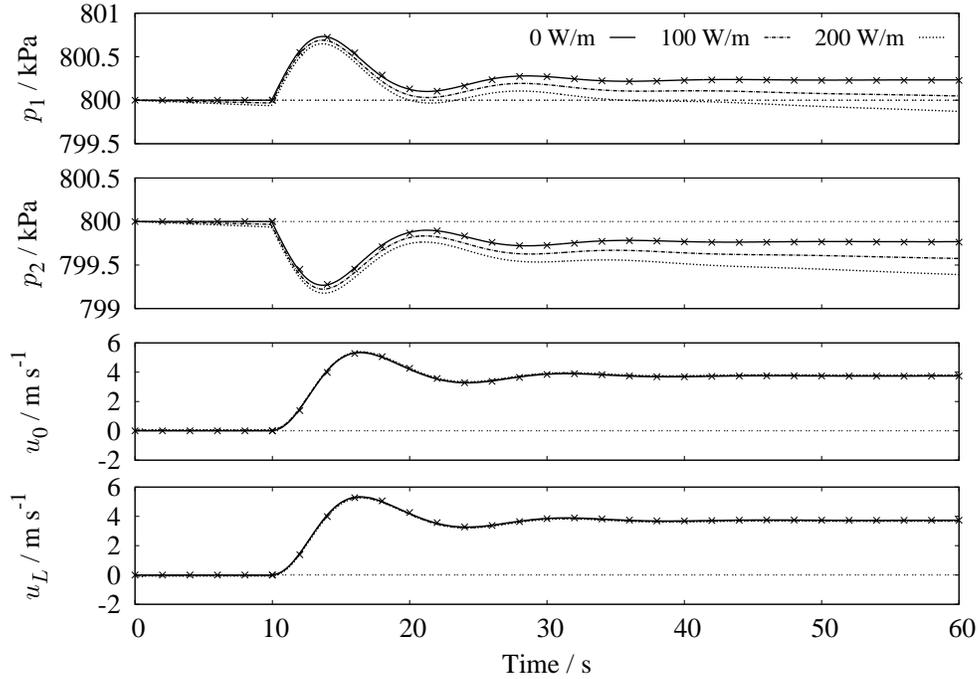}
 \caption{%
 Responses of $p_1$, $p_2$, $u_0:=u(t,0)$, and $u_L:=u(t,L)$ by the original coupled equations. 
 The points ($\times$) show the sequence  of time responses of the derived lumped-parameter model (\ref{eq:network_model}).%
 }
 \label{fig:heatloss}
\end{figure}

Lastly, in order to verify the correctness of the derived model \eref{eq:network_model}, we present a comparison with the brute-force simulation of the original model represented by equations \eref{eq:boiler} to \eref{eq:energy_continuity} and \eqref{eq:load_balance}. 
The equations were implemented using the COMSOL Multiphysics$\textregistered$ Software. 
According to \cite{liu09}, the thermodynamic properties were given by the first order approximation around nominal values. 
From the procedure of lumped-parameter modeling, the approximations needed for deriving \eref{eq:network_model} are the incompressibility condition and the evaluation of thermodynamic quantities in \eqref{eq:thermodynamic_quantities}. 
Since the relevance of these approximations is related to the parameter $\epsilon_3:=d_{\rm r}^2L_{\rm r}Q_{\rm wr}/Q'_{\rm r}$, we present in \figref{fig:heatloss} the simulation results for various setting of the heat loss, represented by $Q_{\rm w}$ in \eqref{eq:energy_continuity}. 
The 
lines in the figure show the results of brute-force simulation of the original model with  $-Q_{\rm w} \cdot \pi d^2/4=0\,{\rm W/m}$, $100\,{\rm W/m}$, and $200\,{\rm W/m}$, respectively. 
The sequence of \emph{points}, denoted by $\times$, represents a sample plot of time responses of \eref{eq:network_model} presented in \figref{fig:short-term}. 
The simulation result clearly shows that the lumped-parameter model \eref{eq:network_model} well approximates the dynamics of the original model when the heat loss is sufficiently small, i.e. when the steam pipes are well insulated. 

For the cases with large heat loss, while 
the pressures $p_1$ and $p_2$ in the original model slowly decrease, the responses of $u_0$ and $u_L$ are well captured by the derived model \eref{eq:network_model}. %
Also, the responses of pressures and velocities in $[10\,{\rm s}, 20\,{\rm s}]$ are correctly produced 
by the derived model. 
The result shows that the derived model describes the transport of steam in the shot-term regime even if the heat loss becomes large. 
This model is thus useful for designing the coordinated operation of heat and electricity supply because the typical time scale of the swing dynamics in power systems is also in the same time regime \cite{machowski08}.

\section{Conclusions}
\label{sec:conclusion}

In this paper, we studied a problem on mathematical modeling for the dynamics in steam supply systems.
The dynamics of interest were originally described by a distributed-parameter model for fast steam flows over a pipe network coupled with a lumped-parameter model for slow internal dynamics of boilers. 
Through physically-relevant approximations, we newly derived a lumped-parameter model that captured stability and multiscale properties of the dynamics.
In order to describe the slow-fast dynamics, we used the notion of Normally Hyperbolic Invariant Manifold (NHIM).
By theoretically analyzing the inner limit of the derived model, we located a set of non-isolated equilibrium points that formed a NHIM. 
Also, the numerical simulations under practical settings of parameters demonstrated the slow-fast dynamics near the NHIM.  
The existence of NHIM clearly suggests the so-called separation principle for operational design of steam supply, which is analogue to power system operation. 

Future directions of this work are as follows. 
One is to verify the correctness and application limit of the derived model via experimental measurements. 
This is inevitable for practical use of the derived model.  
Another one is operational design of multiscale dynamics of steam supply based on the NHIM. 
The characterization is expected to be utilized for separation of different scales or model-order reduction.


\appendix

\section{Derivation of dimensionless governing equations}
\label{sec:scaling}

This first appendix derives the dimensionless equations from \eqref{eq:boiler} to \eref{eq:energy_continuity}. 
The physical quantities with dimension are denoted by superscript $\ast$, and the reference quantities by subscript ${\rm r}$: A physical quantity $z^\ast$ is associated by $z^\ast=z_{\rm r}^\ast z$ to a dimensionless quantity $z$.
Similarly, a function $f^\ast$ is related by $f_{\rm r}^\ast f(z^\ast/z_{\rm r}^\ast)= f^\ast(z^\ast)$ to a dimensionless function $f$. 
The reference quantities are related by the following equations: 
\begin{align}
 &x_{\rm r}^\ast=L_{\rm r}^\ast, \quad  t_{\rm r}^\ast=\dfrac{L_{\rm r}^\ast}{u_{\rm r}^\ast},\quad \rho_{\rm r}^\ast=\rho_{\rm sr}^\ast,\quad h_{\rm r}^\ast=h_{\rm sr}^\ast=h_{\rm wr}^\ast=\dfrac{p_{\rm r}^\ast}{\rho_{\rm r}^\ast},~ \notag \\
 &d_{\rm r}^\ast=\lambda_{\rm r} L_{\rm r}^\ast, \quad m_{\rm r}^{\prime\ast}= \rho_{\rm r}^\ast{d_{\rm r}^\ast}^2 u_{\rm r}^\ast,\quad Q_{\rm r}^{\prime\ast}=Q_{\rm Lr}^{\prime\ast}=h_{\rm r}^\ast\rho_{\rm r}^\ast{d_{\rm r}^\ast}^2 u_{\rm r}^\ast. 
\label{eq:reference_relation}
\end{align}

\subsection{Boiler model}

It is well-known in \cite{astrom00,kim00} that the dynamic behavior of boiler's pressure is well captured by global mass and energy balances. 
The global mass balance is given by 
\begin{align}
 \label{eq:boiler_mass_balance}
  \frac{\rm d}{{\rm d} t^\ast}(\rho_{{\rm s}}^\ast V_{\rm s}^\ast+ \rho_{{\rm w}}^\ast V_{{\rm w}}^\ast) = m_{{\rm f}}^{\prime \ast}-m_{{\rm s}}^{\prime \ast},
\end{align}
and the global mass balance by
\begin{align}
 \label{eq:boiler_energy_balance}
  \frac{\rm d}{{\rm d} t^\ast}\{ (\rho_{{\rm s}}^\ast h_{{\rm s}v}^\ast -p^\ast)V_{{\rm s}}^\ast +(\rho_{{\rm w}}^\ast h_{{\rm w}}^\ast-p^\ast)V_{{\rm w}}^\ast +m_{{\rm t}}^\ast C_{{\rm p}}^\ast T_{{\rm m}}^\ast \} 
 = Q^{\prime\ast} +m_{{\rm f}}^{\prime\ast} h_{{\rm f}}^{\prime\ast} -m_{{\rm s}}^{\prime\ast} h_{{\rm s}}^{\prime\ast},
\end{align}
where the term $h^\ast-p^\ast/\rho^\ast$ corresponds to internal energy. 
Under \aref{as:water_level}, by multiplying \eref{eq:boiler_mass_balance} by $h_{\rm w}^\ast$ and subtracting the result from \eref{eq:boiler_energy_balance} we have
\begin{equation}
 e^\ast \frac{{\rm d} p^\ast}{{\rm d} t^\ast}=Q^{\prime\ast} -m_{\rm f}^{\prime\ast}(h_{\rm w}^\ast-h_{\rm f}^\ast) -m_{\rm s}^{\prime\ast}(h_{\rm s}^\ast-h_{\rm w}^\ast),
  \label{eq:boiler_pressure_dynamics}
\end{equation}
with
\begin{align}
 e^\ast= (h_{\rm s}^\ast-h_{\rm w}^\ast)V_{\rm s}^\ast\frac{\partial\rho_{\rm s}^\ast}{\partial p^\ast} +\rho_{\rm s}^\ast V_{\rm s}^\ast \frac{\partial h_{\rm s}^\ast}{\partial p^\ast} +\rho_{\rm w}^\ast V_{\rm w}^\ast \frac{\partial h_{\rm w}^\ast}{\partial p^\ast} +m_{\rm t}^\ast C_{\rm p}^\ast \frac{\partial T_{\rm m}^\ast}{\partial p^\ast} -V_{\rm s}^\ast -V_{\rm w}^\ast. 
  \label{eq:boiler_coefficient_e}
\end{align}
In addition, \aref{as:boiler} implies $T_{\rm m}^\ast=T_{\rm s}^\ast$, and \aref{as:feedwater} does $h_{\rm f}^\ast=h_{\rm w}^\ast$ in \eref{eq:boiler_pressure_dynamics} and \eref{eq:boiler_coefficient_e}. 
As a result, by using the relation \eref{eq:reference_relation}, the pressure dynamics of boiler are formulated as follows:
\begin{equation}
 e(p) \frac{{\rm d} p}{{\rm d} t}= \dfrac{{d_{\rm r}^\ast}^2 L_{\rm r}^\ast}{e_{\rm r}^\ast} \left\{ Q' -m'_{{\rm s}}h_{{\rm c}}(p) \right\},  \label{eq:boiler_append}
\end{equation} 
Thus, \eref{eq:boiler} is derived by defining the small parameter $\epsilon_1:={d_{\rm r}^\ast}^2 L_{\rm r}^\ast/e_{\rm r}^\ast$.

\subsection{Steam pipe model}

The continuity equations of mass, momentum, and energy with dimension are given as follows \cite{osiadacz87,alobaid08,liu09}: 
\begin{align}
 & \dfrac{\partial\rho^\ast}{\partial t^\ast}+\dfrac{\partial}{\partial x^\ast}(\rho^\ast u^\ast)=0,  \\
 & \dfrac{\partial}{\partial t^\ast}(\rho^\ast u^\ast) + \dfrac{\partial}{\partial x^\ast}(\rho^\ast {u^\ast}^2)
 +\dfrac{\partial p^\ast}{\partial x^\ast} +F_{\rm w}^\ast=0, \\
 & \dfrac{\partial}{\partial t^\ast}\left\{ \rho^\ast \left(h^\ast-\dfrac{p^\ast}{\rho^\ast} \right)\right\} + \dfrac{\partial}{\partial x^\ast}(\rho^\ast u^\ast h^\ast) + Q_{\rm l}^\ast =0, 
\end{align}
where $F_{\rm w}$ stands for the shear force acting on a steam element and is approximated by the Darcy-Weisbach equation \cite{osiadacz87,bergman11,kanro:english} as follows:
\begin{equation}
 F_{\rm w}^\ast := \lambda\dfrac{\rho^\ast u^\ast|u^\ast|}{2d^\ast}, \label{eq:darcy-weisbach}
\end{equation}
The coefficient $\lambda$ depends on the Reynolds number $Re$, pipe diameter $d$, and roughness of the inner surface of the pipe and therefore varies according to the steam velocity $u$.
In the case of laminar flow under a low Reynolds number, the value of $\lambda$ is derived in \cite{bergman11} theoretically as $\lambda=64/Re$.
On the other hand, in the case of turbulent flow under high Reynolds numbers, $\lambda$ is approximated by a constant that is determined by $d$ and the roughness \cite{kanro:english,bergman11}. 
In this paper, the friction coefficient $\lambda$ is considered as a constant because the steam flow used in standard steam supply is turbulent \cite{kanro:english,bergman11}.
By using the relation \eref{eq:reference_relation}, the above equations become
\begin{align}
 & \dfrac{\partial\rho}{\partial t}+\dfrac{\partial}{\partial x}(\rho u)=0, \\
 & \dfrac{\partial}{\partial t}(\rho u) + \dfrac{\partial}{\partial x}(\rho u^2)
 +\dfrac{p_{\rm r}^\ast}{\rho_{\rm r}^\ast{u_{\rm r}^\ast}^2} \dfrac{\partial p}{\partial x} +\lambda \dfrac{\rho_l u|u|}{2d} =0, \\
 & \dfrac{\partial}{\partial t}(\rho h)+ \dfrac{\partial}{\partial x}(\rho h u)  = \dfrac{\partial p}{\partial t} + \dfrac{ {d_{\rm r}^\ast}^2 L_{\rm r}^\ast Q_{\rm wr}^\ast}{Q_{\rm r}^{\prime\ast}} Q_{{\rm w}}. 
\end{align}
Thus, \eref{eq:mass_continuity} to \eref{eq:energy_continuity} are derived by defining the parameters $\epsilon_2:=\rho_{\rm sr}^\ast{u_{\rm r}^\ast}^2/p_{\rm r}^\ast$ and $\epsilon_3:={d_{\rm r}^\ast}^2L_{\rm r}^\ast Q_{\rm wr}^\ast/Q_{\rm r}^{\prime\ast}$.

\section{Summarized Graph Theory}
\label{sec:gt}
The second appendix provides a summarized theory of graph from \cite{iri69}.  
Along the notation introduced in \secref{ssec:notation}, 
consider a directed graph $\mathcal{G}$ with $n$ vertices and $m$ links. 
Assume $\mathcal{G}$ is connected and is represented by the incidence matrix ${\sf R}$.  
The matrix %
 ${\sf R}$ (or ${\sf R}^{\top}$) is regarded as a linear map ${\sf R}: \mathbb{R}^m \to \mathbb{R}^n$ (or ${\sf R}^\top: \mathbb{R}^n \to \mathbb{R}^m$), and its image and kernel are related to the graph's topology. 
$\Ker({\sf R})$ and ${\rm Im}({\sf R}^\top)$ are subspaces of $\mathbb{R}^m$, and their dimensions coincide with the number of independent loops and cutsets, respectively. 
Since $\mathcal{G}$ is connected, we have ${\rm dim} (\Ker({\sf R})) = m-n+1$ and ${\rm dim (Im}({\sf R}^\top))= n-1$ (see \cite{iri69}), where ${\rm dim}(X)$ stands for the dimension of $X$. 
Also, ${\rm Ker}({\sf R}^\top)$ is a subspace of $\mathbb{R}^n$ 
given as 
\begin{equation}
  {\rm \Ker}({\sf R}^\top) = \{ c {\bm 1} \in \mathbb{R}^n ~ |~ c \in \mathbb{R} \}.
 \label{eq:kerRt}
\end{equation}
This result is derived from  ${\rm dim (Im}({\sf R}^\top))= n-1$ and the fact that every link connects exactly two vertices.
Finally, ${\rm Im}({\sf R)}$ is a subspace of $\mathbb{R}^n$ orthogonal to ${\rm Ker}({\sf R}^\top)$, and its dimension is equal to $n-1$.


\listoftables

\listoffigures

\end{document}